\documentclass[12pt, a4paper]{article}

\usepackage{a4wide}
\usepackage{amsmath}
\usepackage{amsfonts}
\usepackage{enumerate}
\date{}
\begin{document}
\title{\large\textbf{ON A THEOREM BY PELC AND PRIKRY ON THE\\ NON EXISTENCE OF INVARIANT EXTENSIONS OF\\ BOREL MEASURES }}
\author{\normalsize
\textbf{S.BASU \& D.SEN}\\
}
\date{}	
\maketitle
{\small \noindent \textbf{{\textbf{AMS subject classification ($2010$):}}}} {\small $28A05$, $28D05$, $28D99$.\\}
{\small\noindent \textbf{\textbf{\textmd{\textbf{{Keywords and phrases :}}}}} Banach-Kuratowski Matrix, diffused admissible functional, space with transformation group, $k$-small system, upper semicontinuous $k$-small system, $k$-additive(weakly $k$-addtive)algebra admissible with respect to a $k$-small system, generalized continuum hypothesis.}\\
\vspace{.03cm}\\
{\normalsize
\textbf{\textbf{ABSTRACT:}} There are certain countably generated $\sigma$-algebras of sets in the real line which do not admit any non-zero, $\sigma$-finite, diffused (or, continuous) measure. Such countably generated $\sigma$-algebras can be obtained by the use of some special types of infinite matrix known as the Banach-Kuratowski matrix and the same may be used in deriving a generalized version of Pelc and Prikry's theorem as shown by Kharazishvili. In this paper, using some methods of combinatorial set theory and some modified version of the notion of small sets originally introduced by Rie$\check{c}$an, Rie$\check{c}$an and Neubrunn, we give an abstract and generalized formulation of Pelc and Prikry's theorem in spaces with transformation groups.\\
\begin{center}
\section{\large{INTRODUCTION}}
\end{center}
\normalsize Banach $[1]$ (see also $[4]$) asked if there exist two countably generated $\sigma$-algebras on the interval $[0,1)$ such that they both carry probability diffused measures, whereas the $\sigma$-algebra generated by their union does not. An answer to this was provided by Grzegorek in $[2]$ using Martin axiom and in $[3]$ (see, also $[4]$) without using any such additional set theoretic assumptions. In $[10]$, Pelc and Prikry obtained an analogue of the result in translation invariant settings, and, Kharazishvili $[7]$ obtained a generalization of Pelc and Prikry's result by constructing (under continuum hypothesis) certain Banach-Kuratowski matrix consisting of sets that are almost invariant with respect to the group of all isometric transformations.\
\vspace{.001cm}

Let $\omega$ and $\omega_{_{1}}$ denote the first infinite and the first uncountable ordinals and $F= \omega^{\omega}$ be the family of all functions from $\omega$ into $\omega$. For any two functions $f$ and $g$ from $F$, we write $f\preceq g$ to mean that there exists a natural number $n(f,g)$ such that $f(m)\leq g(m)$ for all $m$ such that $n(f,g)\leq m$. The relation $\preceq$ so defined is a pre ordering on $F$ and under the assumption of continuum hypothesis, it is not hard to define a subset $E= \{f_{\xi} : \xi< \omega_{_{1}}\}$ of $F$ satisfying the following two conditions:\\
(i) If $f$ is any arbitrary function from $F$, then there exists an ordinal $\xi< \omega_{_{1}}$ such that $f\prec f_{\xi}$. In otherwords, $E$ is cofinal in $F$.\\
(ii) For no two ordinals $\xi$ and $\rho$ satisfying $\xi< \rho< \omega_{_{1}}$ does the relation $f_{\rho}\prec f_{\xi}$ holds true.\\
\vspace{.03cm}

Now conditions (i) and (ii) imply that card${(E)}= \omega_{_{1}}$. Further, if for any two natural numbers $m$ and $n$, we set $E_{m}{_,}{_{n}}= \{f_{\xi}\in E : f_{\xi}(m)\leq n\}$, then the countable double family of sets $(E_{m}{_,}{_{n}})_{m<\omega}{_,}{_{n<\omega}}$ satisfies\\
(a) $E_{m}{_,}{_{0}}\subseteq E_{m}{_,}{_{1}}\subseteq \ldots\ldots\subseteq E_{m}{_,}{_{n}}\subseteq \ldots\ldots$ for any natural number $m<\omega$.\\
(b) $E= \cup \{E_{m}{_,}{_{n}} : n< \omega\}$ and\\
(c) $E_{0}{_,}{_{f(0)}}\cap E_{1}{_,}{_{f(1)}}\cap \ldots\ldots E_{m}{_,}{_{f(m)}}\cap\ldots$ is atmost countable for every $f\in F$.\\
\vspace{.001cm}

A matrix $(E_{m}{_,}{_{n}})_{m<\omega}{_,}{_{n<\omega}}$ on $E$ having the above three properties is called a Banach-Kuratowski matrix $[5]$, and it can be proved that there does not exist any non-zero, $\sigma$-finite, diffused (or, continuous) measure defined simultaneously for all sets $E_{m}{_,}{_{n}}$. Not only that, the existence of a Banach-Kuratowski matrix on $E$ proves even more $[5]$ and it is this that there does not exist any non-zero, diffused, admissible functional defined simultaneously for all the above sets $E_{m}{_,}{_{n}}$, where by a diffused admissible functional $[5]$ (see also $[7]$) we mean a set valued mapping $\nu$ defined on a family of subsets of $E$ which is closed under finite intersection and for which the following set of conditions are fulfilled :\\
(1) $\nu$ is defined on every countable set $X\subseteq E$ with $\nu(X)= 0$.\\
(2) If $\{Z_{n}: n< \omega\}$ is an increasing family of sets (with respect to inclusion) from the domain of $\nu$, then their union $\displaystyle{\bigcup_{n<\omega}}{Z_{n}}$ is also a member of the domain and \begin{center} $\nu(\displaystyle{\bigcup_{n<\omega}}{Z_{n}}) \leq sup\{\nu(Z_{n}): n<\omega\}$
\end{center}
(3) If $\{Z_{n}: n< \omega\}$ is a decreasing family of sets (with respect to inclusion) from the domain of $\nu$, then their intersection $\displaystyle{\bigcap_{n<\omega}}{Z_{n}}$ is also a member of the domain and \begin{center} $\nu(\displaystyle{\bigcap_{n<\omega}}{Z_{n}}) \geq inf\{\nu(Z_{n}): n<\omega\}$
\end{center}
Evidently, every finite measure defined on a $\sigma$-algebra of sets satisfies conditions $(1)-(3)$ whereas an admissible functional is not necessarily a measure for it need not have the property of $\sigma$-additivity as required in the definition of any ordinary measure function.\\
\vspace{.001cm}

Pelc and Prikry $[10]$ under the assumption of continuum hypothesis and the use of Hamel basis, proved that\\
\vspace{.01cm}

\textbf{THEOREM PP :} \textit{There exist countably generated $\sigma$-algebras $\mathfrak R_{1}$, $\mathfrak R_{2}$ of subsets of $[0,1)$ and probability measures $\mu_{_{1}}$ and $\mu_{_{2}}$ on $\mathfrak R_{1}$, $\mathfrak R_{2}$ respectively such that \\
(i) $\mathfrak R_{1}$, $\mathfrak R_{2}$ both contain all Borel sets and are translation invariant. \\
(ii) $\mu_{_{1}}$ and $\mu_{_{2}}$ both extend the Lebesgue measure and are translation invariant.\\
(iii) there is no non-atomic probability measure on any $\sigma$-algebra containing $\mathfrak R_{1}\cup \mathfrak R_{2}$.\\
(Here non-atomic measure means that it is diffused or continuous) } \
\vspace{.001cm}

Kharazishvili $[7]$ by constructing certain Banach-Kuratowski matrix consisting of almost invariant sets (with respect to the group $\Gamma$ of all isometric transformations of $\mathbb{R}$) obtained the following generalization of Pelc and Prikry's theorem.\\
\vspace{.01cm}

\textbf{THEOREM K :}\textit{ Suppose the continuum hypothesis hold. Then there are two $\sigma$-algebras $\mathcal S_{1}$ and $\mathcal S_{2}$ of sets in $\mathbb{R}$ having the following properties : \\
(i) $\mathcal S_{1}$ and $\mathcal S_{2}$ are countably generated and invariant under $\Gamma$.\\
(ii) the $\sigma$-algebra of Borel sets $\mathcal{B}(\mathbb{R})\subseteq \mathcal S_{1}\cap \mathcal S_{2}$\\
(iii) there exists a $\Gamma$-invariant measure $\mu_{_{1}}$ on $\mathcal S_{1}$ extending the standard Borel measure on $\mathbb{R}$.\\
(iv) there exists a $\Gamma$-invariant measure $\mu_{_{2}}$ on $\mathcal S_{2}$ extending the standard Borel measure on $\mathbb{R}$.\\
(v) there is no non-zero, diffused admissible functional defined on the $\sigma$-algebra generated by $\mathcal S_{1}\cup \mathcal S_{2}$.}\
\vspace{.01cm}

In this paper, our aim is to give an abstract version of Pelc and Prikry's theorem in spaces with transformation groups. Instead of using measure or admissible functional, we utilize some more general concepts which are introduced in the following section. \\
\newpage
\begin{center}
\section{\large{PRELIMINARIES AND RESULTS}}
\end{center}
\normalsize A space $X$ with a transformation group $G$ is a pair $(X,G)$ $[4]$ where $X$ is a nonempty basic set and $G$ is a group acting on $X$ which means that there is a function $(g, x)\rightarrow gx$ from $G\times X$ onto $X$ satisfying the following two conditions:\\
TG 1) for each $g\in G$, $x\rightarrow gx$ is a bijection (or, permutation) of $X$\\
TG 2) for all $x\in X$, and $g_{_{1}}, g_{_{2}}\in G$, $g_{_{1}}(g_{_{2}}x)= g_{_{1}}g_{_{2}}x$.\\
\vspace{.01cm}

We say that $G$ acts freely on $X$ $[4]$ if $\{x \in X : gx = x\} = \emptyset$ for all $g\in G\setminus \{e\}$ where $e$ is the identity element of $G$ (in fact, the identity element $`e'$ can be considered as the identity transformation $e : X\rightarrow X$ on X). For any $g \in G$ and $E \subseteq X$, we write $gE$ $[4]$ to denote the set $\{gx : x \in E\}$ and call a nonempty family (or, class) $\mathcal A$ of sets $G$-invariant $[4]$ if $gE\in\mathcal A$ for every $g\in G$ and $E\in\mathcal A$. If $\mathcal A$ is a $\sigma$-algebra, then a measure $\mu$ on $\mathcal A$ is called $G$-invariant $[4]$ if $\mathcal A$ is a $G$-invariant class and $\mu(gE) = \mu(E)$ for every $g\in G$ and $E\in \mathcal A$. It is called $G$-quasiinvariant $[4]$ if $\mathcal A$ and the $\sigma$-ideal generated by $\mu$-null sets are both $G$-invariant classes. Obviously, any $G$-invariant measure is also $G$-quasiinvariant but not conversely. Any set of the form $Gx = \{gx : g\in G\}$ for some $x\in X$ is called a $G$-orbit $[4]$  of $x$. The collection of all $G$-orbits give rise to a partition of $X$ into mutually disjoint sets. A subset $E$ of $X$ is called a complete $G$-selector $[4]$ (or, simply, a $G$-selector) in $X$ if $E\hspace{.0005cm}\cap\hspace{.0005cm}Gx$ consists of exactly one point for each $x\in X$. If $Gx= E$ for each $x\in X$, then $G$ is said to act transitively on $X$ $[4]$. In this situation, for any $x,y\in X$, there exists $g\in G$ such that $y= gx$.\
\vspace{.01cm}

Throughout this paper, we identify every infinite cardinal with the least ordinal representing it and every ordinal with the set of all ordinals preceeding it. We write card$(A)$ and card$(\mathcal A)$ to represent the cardinality of any set $A$ or any family (or, class) of sets $\mathcal A$ and use symbols such as $\alpha$, $\beta$, $\gamma$, $\delta$, $\xi$, $\eta$, $\rho$, $k$ etc for infinite cardinals. Moreover, we denote by $k^{^{+}}$ the successor of $k$. Now given any space $(X,G)$ with a transformation group $G$ and an infinite limit cardinal $k$, we define\\
\vspace{.01cm}

\textbf{DEFINITION $\textbf{2.1}$ :} A nonempty class $\mathcal S$ of subsets of $X$ as a $k$-additive algebra on $(X,G)$ if\\
(i) $\mathcal S$ is an algebra\\
(ii) $\mathcal S$ is $k$-additive which means that $\mathcal S$ is closed with respect to the union of atmost $k$ number of sets from it.\\
(iii) $\mathcal S$ is a $G$-invariant class.\
\vspace{.01cm}

A nonempty class $\mathcal S$ is called weakly $k$-additive if in statement (ii) above, we replace `` atmost $k$ '' by `` by less than $k$ '' number of sets.\
\vspace{.01cm}

 Thus a $k$-additive (or, a weakly $k$additive) algebra on $(X,G)$ is a $k$-additive (or, a weakly $k$-additive) algebra on $X$ which is also $G$-invariant. A $k$-additive algebra (or, a weakly $k$-additive) algebra on $X$ is called diffused if every singleton set $\{x\}\in \mathcal S$.\
 \vspace{.01cm}

 The notion of ``small system'' or ``system of small sets'' was originally introduced by Rie$\acute{c}$an, Rie$\acute{c}$an and Neubrunn $[11]$, $[12]$, $[14]$ and used by many authors $[9]$, $[13]$ etc to establish abstract formulations of several well known classical theorems of Lebesgue measure and integration. For our purpose, we need to have a modified version of this concept.\
 \vspace {.01cm}

 By a $k$-small system on $(X,G)$,we mean\\

\textbf{DEFINITION $\textbf{2.2}$ :} A transfinite $k$-sequence $\{\mathcal N_{_{\alpha}}\}_{_{\alpha<k}}$ where each $\mathcal N_{_{\alpha}}$ is a class of subsets of $X$ satisfying the following set of conditions :\\
(i) $\emptyset\in \mathcal N_{\alpha}$ for all $\alpha<k$.\\
(ii) Each $\mathcal N_{\alpha}$ is a $G$-invariant class.\\
(iii) $E\in \mathcal N_{\alpha}$ and $F\subseteq E$ implies $F\in \mathcal N_{\alpha}$. That is, $\mathcal N_{_{\alpha}}$ is a hereditary class.\\
(iv) $E\in \mathcal N_{\alpha}$ and $F\in{\displaystyle{\bigcap_{\alpha<k}}\hspace{.01cm}{\mathcal N_{\alpha}}}$ implies $E\cup F\in\mathcal N_{\alpha}$\\
(v) For any $\alpha < k$, there exists $\alpha^{\ast}> \alpha$ such that for any one-to-one correspondence $\beta\rightarrow \mathcal N_{_{\beta}}$ with $\beta > \alpha^{\ast}$, ${\displaystyle{\bigcup_{\beta}{E_{_{\beta}}}}}\in \mathcal N_{_{\alpha}}$ whenever $E_{_{\beta}}\in \mathcal N_{_{\beta}}$.\\
(vi)  For any $\alpha , \beta < k$, there exists $\gamma > \alpha , \beta$  such that $\mathcal N_{_{\gamma}} \subseteq \mathcal N_{ _{\alpha}}$ and  $\mathcal N _{_{\gamma}} \subseteq \mathcal N _{_{\beta}}$; i.e $\mathcal N_{_{\alpha}}$ is directed.\
\vspace{.01cm}\\

\textbf{DEFINITION $\textbf{2.3}$ :} A $k$-additive (or,weakly $k$-additive) algebra $\mathcal S$ is admissible with respect to a $k$-additive small system $\{\mathcal N_{_{\alpha}}\}_{_{\alpha<k}}$ if\\
(i) $\mathcal S\setminus\mathcal N_{\alpha}\neq\emptyset$ for some $\alpha<k$ and $\mathcal S\cap\mathcal N_{\alpha}\neq \emptyset$ for every $\alpha<k$. In otherwords, $\mathcal S$ is compatible with respect to $\{\mathcal N_{_{\alpha}}\}_{_{\alpha<k}}$.\\
(ii) $\mathcal N_{\alpha}$ has a $\mathcal S$-base i.e $E\in\mathcal N_{\alpha}$ is contained in some $F\in \mathcal N_{\alpha}\cap\mathcal S$,\\
and (iii) $\mathcal S\setminus\mathcal N_{\alpha}$ satisfies $k$-chain condition, i.e, the cardinality of any arbitrary collection of mutually disjoint sets from $\mathcal S\setminus\mathcal N_{\alpha}$ is atmost $k$.\
\vspace{.01cm}

 By virtue of conditions (i)-(iii) and (v) of Definition $2.2$, it follows that $\mathcal N_{\infty}= {\displaystyle{\bigcap_{\alpha<k}}}\hspace{.02cm}{\mathcal N_{\alpha}}$ is a $k$-additive ideal on $(X,G)$ which means that the set $\mathcal N_{\infty}$ is a $k$-additive ideal in $X$ and also a $G$-invariant class.\
\vspace{.01cm}

To the list of the above definitions, we further add that\\
\vspace{.01cm}

\textbf{DEFINITION $\textbf{2.4}$ :} A $k$-small system $\{\mathcal N_{_{\alpha}}\}_{_{\alpha<k}}$ is upper semicontinuous relative to a $k$-additive (or, a weakly $k$-additive) algebra $\mathcal S$ on $(X,G)$ if for every nested $k$-sequence $\{E_{_{\xi}}:\xi<k\}$ of sets from $\mathcal S$ satisfying $E_{_{\xi}}\notin \mathcal N_{_{\alpha}}{_{_{_{0}}}}$ for some $\alpha_{_{0}}<k$ and all $\xi<k$, we have $\displaystyle{\bigcap_{\xi}}E_{_{\xi}}\notin \mathcal N_{\infty}$.\\
\vspace{.01cm}

\textbf{DEFINITION $\textbf{2.5}$ :} A set $E\subseteq X$ is $(\mathcal S, \{\mathcal N_{_{\alpha}}\}_{_{\alpha<k}})$-thick in $X$ if $B\subseteq {X\setminus E}$ and $B\in\mathcal S$ implies $B\in \mathcal N_{\infty}$. More generally, $E$ is $(\mathcal S, \{\mathcal N_{_{\alpha}}\}_{_{\alpha<k}})$-thick in $F(\supseteq E)\in\mathcal S$ if $B\subseteq {F\setminus E}$ and $B\in\mathcal S$ implies that $B\in\mathcal N_{\infty}$ and $F$ is then called a $(\mathcal S, \{\mathcal N_{_{\alpha}}\}_{_{\alpha<k}})$-cover of $E$ in $X$. If $F_{_{1}}$  and $F_{_{2}}$ are any two $(\mathcal S, \{\mathcal N_{_{\alpha}}\}_{_{\alpha<k}})$-covers of $E$ then their symmetric difference $F_{_{1}}\Delta F_{_{2}}\in \mathcal N_{\infty}$. Thus in a sense they are identical and so from now on, instead of using `` a $\ldots$ -cover" we will use the phrase `` the $\ldots$ -cover" as and when required.\
\vspace{.01cm}

Using admissibility of $\mathcal S$ and the above definitions, we prove the following\\
\vspace{.01cm}

\textbf{PROPOSITION} \textbf{$\textbf{2.6}$ :} \textit{If $E\notin\mathcal N_{\infty}$ and $gE\Delta E\in\mathcal N_{\infty}$ for every $g\in G$, then $E$ is $(\mathcal S, \{\mathcal N_{_{\alpha}}\}_{_{\alpha<k}})$-thick in $X$.
}\\
\vspace{.1cm}

\textbf{PROOF :} Let $B\subseteq {X\setminus E}$ and suppose $B\in {\mathcal S\setminus \mathcal N_{\infty}}$. Then because of admissibility of $\mathcal S$ with respect to $\{\mathcal N_{_{\alpha}}\}_{_{\alpha<k}}$, it is possible to generate (using $k$-chain condition) a $k$-sequence $\{g_{_{\alpha}}:\alpha<k\}$ in $G$ such that $X\setminus {\displaystyle\bigcup_{\alpha<k}}{g_{_{\alpha}}B}\in \mathcal N_{\infty}$. But then there exists $\alpha_{_{0}}<k$ such that ${g_{\alpha_{_{0}}}B\cap E\notin \mathcal N_{\infty}}$ for otherwise $E\in \mathcal N_{\infty}$ because $\mathcal N_{\infty}$ is $k$-additive (we have already noted that $\mathcal N_{\infty}$ is a $k$-additive ideal on $(X,G)$). Hence $B\cap {{g^{-1}_{_{\alpha_{_{0}}}}}}E\notin \mathcal N_{\infty}$ by condition (ii) of Definition $2.2$. But this contradicts the hypothesis. Therefore $B\in \mathcal N_{\infty}$.\\
\vspace{.1cm}

In the family $F$ of all functions from $k$ into $k$, let us set up a preordering as follows:\\
$f,g\in k^{k}$, $f\preceq g$ iff card$(\{\alpha< k : g(\alpha)<f(\alpha)\})<k$. If the generalized continuum hypothesis is assumed, it is not hard to define a subfamily $E= \{f_{\xi} : \xi< k^{^{+}}\}$ of $F$ which is cofinal in $F$ in the sense that for every $f\in F$, there is some $f_{\xi}\in E$ such that $f\preceq f_{\xi}$ and which also satisfies the property that the relation $f_{\eta}\prec f_{\xi}$ is not true for any $\xi$ and $\eta$ such that $\xi< \eta$.\\
Now upon setting $E_{\alpha}{_,}{_{\beta}}=\{f_{\xi}\in E : f_{\xi}(\alpha)\leq \beta\}$, we find that the double family $(E_{\alpha}{_,}{_{\beta}})_{_{\alpha< k}}{_,}{_{_{\beta< k}}}$ of sets satisfies the following three conditions :\\
(i) $E_{\alpha}{_,}{_{\beta}}\subseteq E_{\alpha}{_,}{_{\gamma}}$ for any $\alpha< k$ and $\beta\leq\gamma< k$.\\
(ii) $E=\displaystyle{\bigcup_{\beta< k}{{E_{\alpha}{_,}{_{\beta}}}}}$ for any $\alpha<k$, and\\
\vspace{.001cm}\\
(iii) for any $f\in k^{^{k}}$, card $(\displaystyle{\bigcap_{\alpha<k}\hspace{.0001cm}{{{E_{\alpha}{_,}{_{f(\alpha)}}}}}})\leq k$.\
\vspace{.01cm}

It may be noted that the family $(E_{\alpha}{_,}{_{\beta}})_{_{\alpha< k}}{_,}{_{_{\beta< k}}}$ constructed above is a direct generalization of the Banach-Kuratowski matrix. In the following two theorems, we assume generalized continuum hypothesis.\\
\vspace{.01cm}

\textbf{THEOREM $\textbf{2.7}$ :} \textit{Let $(X,G)$ be a space with a transformation group $G$ where $k^{^{+}}$= card $G$ $\leq$ card $X$ and $G$ acts freely on $X$. Let $L$ be a $G$-selector in $X$. Then there exists a family $(F_{\alpha}{_,}{_{\beta}})_{_{\alpha< k}}{_,}{_{_{\beta< k}}}$ of sets in $X$ which is not contained in any $k$-additive algebra $\mathcal S$ on $(X,G)$ which contains $L$ and which is admissible with respect to any $k$-small system $\{\mathcal N_{_{\alpha}}\}_{_{\alpha<k}}$ on $(X,G)$ which is upper semicontinuous relative to $\mathcal S$.}\\
\vspace{.1cm}

\textbf{PROOF :} We write $G=\displaystyle{\bigcup_{\rho< k^{^{+}}}{G_{\rho}}}$ where $\{G_{\rho} : \rho< k^{^{+}}\}$ is an increasing family of subgroups of $G$ satisfying $G_{\rho}\neq\displaystyle{\bigcup_{\eta< \rho}}{G_{\eta}}$ and card{{($G_{\rho}$)}}$\leq$ $k$ for every $\rho< k^{^{+}}$ (for this representation, see $[6]$, Exercise $19$, Ch $3$). Since by hypothesis, $G$ acts freely on $X$ so the above increasing family yields a disjoint covering $\{\Omega_{\gamma} : \gamma< k^{^{+}}\}$ of $X$ where $\Omega_{\gamma}= (G_{\gamma}\setminus\displaystyle{\bigcup_{\eta< \gamma}}{G_{\eta}})L$. Using the $k\times k$ matrix $(E_{\alpha}{_,}{_{\beta}})_{_{\alpha< k}}{_,}{_{_{\beta< k}}}$, we define ${F_{\alpha}{_,}{_{\beta}}=\displaystyle{\bigcup_{\gamma\in E_{\alpha}{_,}{_{\beta}}}}{\Omega_{\gamma}}}$. Then $F_{\alpha}{_,}{_{\beta}}\subseteq F_{\alpha}{_,}{_{\gamma}}$ whenever $\beta< \gamma< k$ and $X=\displaystyle{\bigcup_{\beta< k}}{F_{\alpha}{_,}{_{\beta}}}$.\
\vspace{.1cm}

We observe that there exists $\delta$ such that no subset $M$ of $X$ can belong to $\mathcal N_{\delta}$ if its complement in $X$ i.e $X\setminus M$ is in $\mathcal N_{\delta}$. This follows since $\mathcal S\setminus\mathcal N_{_{\alpha}}\neq\emptyset$ for some $\alpha<k$ (by (i) of Definition $2.3$) and so $X\in\mathcal S\setminus\mathcal N_{_{\alpha}}$ (by (iii) of Definition $2.2$) and also because there exists $\beta, \gamma > \alpha$ (by (v) of Definition $2.2$) such that $\mathcal N_{\beta}\cup\mathcal N_{\gamma}\subseteq\mathcal N_{\alpha}$ and $\delta>\beta,\gamma$ (by (vi) of Definition $2.2$) such that $\mathcal N_{\delta}\subseteq \mathcal N_{\beta}$, $\mathcal N_{\gamma}$. Now, suppose $\mathcal S$ be any $k$-additive algebra on $(X,G)$ which is admissible with respect to a $k$-small system $\{\mathcal N_{_{\alpha}}\}_{_{\alpha<k}}$ on $(X,G)$ which is upper semicontinuous relative to $\mathcal S$, and such that $L\in \mathcal S$.\\
If possible, let $\{F_{\alpha}{_,}{_{\beta}} : \alpha<k, \beta<k\}\subseteq \mathcal S$.
\vspace{.1cm}

By virtue of condition (v) of Definition $2.2$ and also because $\{\mathcal N_{_{\alpha}}\}_{_{\alpha<k}}$ is upper semicontinuous relative to $\mathcal S$, there exists $\delta^{^{\ast}}>\alpha$ and an one-to-one correspondence $f: \alpha\rightarrow\beta_{_{\alpha}}$ $(\beta_{_{\alpha}}>\delta^{^{\ast}})$ such that ${G_{\alpha}{_,}{_{\beta_{_{\alpha}}}}}\in {\mathcal N_{\beta_{_{\alpha}}}\cap \mathcal S}$ and $\displaystyle{\bigcup_{\alpha< k}}{G_{\alpha}{_,}{_{\beta_{_{\alpha}}}}}\in{\mathcal N_{\delta}\cap\mathcal S}$ where ${G_{\alpha}{_,}{_{\beta_{_{\alpha}}}}}= X\setminus {{F_{\alpha}{_,}{_{\beta_{_{\alpha}}}}}}$. Therefore, $\displaystyle{\bigcap_{\alpha< k}}{F_{\alpha}{_,}{_{\beta_{_{\alpha}}}}}\in{\mathcal S\setminus\mathcal N_{_{\delta}}}$ by virtue of the above observation. But $\displaystyle{\bigcap_{\alpha< k}}{F_{\alpha}{_,}{_{f(\alpha)}}}$ is the union of atmost $k$ translates of $L$. But $L$ cannot belong to the complement of any $\mathcal N_{\alpha}$ because card($G$)= $k^{^{+}}$, $G$ acts freely on $X$ and $\mathcal S$ being admissible with respect to $\{\mathcal N_{_{\alpha}}\}_{_{\alpha<k}}$ satisfies the $k$-chain condition. Hence, $L\in \mathcal S\cap \mathcal N_{\infty}$ and consequently, $\displaystyle{\bigcap_{\alpha< k}}{F_{\alpha}{_,}{_{f(\alpha)}}}\in{\mathcal S\cap \mathcal N_{\infty}}$ since $\mathcal N_{\infty}$ is a $k$-additive ideal on $(X,G)$. Thus $\displaystyle{\bigcap_{\alpha< k}}{F_{\alpha}{_,}{_{f(\alpha)}}}\in{\mathcal S\setminus{\mathcal N_{_{\delta}}}}$ and also $\displaystyle{\bigcap_{\alpha< k}}{F_{\alpha}{_,}{_{f(\alpha)}}}\in{\mathcal S\cap{\mathcal N_{\infty}}}$ - a contradiction.\\
\vspace{.001cm}

This proves the theorem.\\
\vspace{.01cm}

We now give an abstract formulation of Pelc and Prikry's theorem.\\
\vspace{.1cm}

\textbf{THEOREM $\textbf{2.8}$ :} \textit{Let $(X,G)$ be a space with transformation group $G$, $k$ be any regular infinite cardinal such that $k^{^{+}}= card(G)= card(X)$. Assume that $G$ acts freely and transitively on $X$. Let $\mathcal S_{0}$ be a diffused $k$-additive algebra on $(X,G)$ admissible with respect to a $k$-small system $\{\mathcal N_{_{\alpha}}{{^{^{0}}}}\}_{_{\alpha<k}}$ on $(X,G)$ which is upper semicontinuous relative to $\mathcal S_{0}$ and such that $\mathcal N_{_{\alpha}}{{^{^{0}}}}\subseteq\mathcal S_{0}$ for every $\alpha<k$. Moreover, let every set in $X$ possesses a $(\mathcal S_{0},\{\mathcal N_{_{\alpha}}{{^{^{0}}}}\}_{_{\alpha<k}})$-cover. Then there exist diffused weakly $k$-additive algebras $\mathcal S_{1}$ and $\mathcal S_{2}$ on $(X,G)$ such that \\
(i) both $\mathcal S_{1}$ and $\mathcal S_{2}$ properly contain $\mathcal S_{0}$ and are $k$-generated (i.e generated by a class of cardinality $k$) if $\mathcal S_{0}$ is so.\\
(ii) $\mathcal S_{1}$ and $\mathcal S_{2}$ are each admissible with respect to some upper semicontinuous $k$-small system on $(X,G)$.\
}

\textit{But no $k$-additive algebra on $(X, G)$ containing $\mathcal S_{1}\cup\mathcal S_{2}$ exists which is admissible with respect to any upper semicontinuous $k$-small system on $(X,G)$.}\\
\vspace{.1cm}

\textbf{PROOF :}  Since $G$ acts freely and transitively on $X$, so for any arbitrary but fixed choice of $x$ from $X$, the increasing family $\{G_{\rho} : \rho< k^{^{+}}\}$ (as constructed in the proof of theorem $2.7$) of subgroups yields a disjoint covering $\{\Lambda_{\gamma} : \gamma< k^{^{+}}\}$ of $X$ where $\Lambda_{\gamma}= (G_{\gamma}\setminus\displaystyle{\bigcup_{\eta< \gamma}}{G_{\eta}})x$. Again as $\{x\}\in \mathcal S_{0}$ and $\mathcal S_{0}$ is admissible with respect to $\{\mathcal N_{_{\alpha}}{{^{^{0}}}}\}_{_{\alpha<k}}$, so by virtue of condition (iii) of Definition $2.3$, $\{x\}\in\mathcal N_{_{\infty}}{^{^{0}}}$. Now by theorem $2.7$ with $L$ replaced by $\{x\}$, we find that $\{\mathcal H_{\alpha}{_,}{_{\beta}} : \alpha<k, \beta<k\}\not\subseteq \mathcal S_{0}$ where $\mathcal H_{\alpha}{_,}{_{\beta}}=\displaystyle{\bigcup_{\gamma\in E_{\alpha}{_,}{_{\beta}}}}{\Lambda_{\gamma}}$.
Consequently, there exists $\mathcal H_{\alpha_{0}}{_,}{_{\beta_{0}}}\notin \mathcal S_{0}$.\\
\vspace{.01cm}

From the constructions of $\Lambda_{\gamma}$, it readily follows that for any $\Sigma\subseteq k^{^{+}}$,\\ card $(g({\displaystyle\bigcup_{\gamma\in\Sigma}}{\Lambda_{\gamma}})\Delta\ ({\displaystyle\bigcup_{\gamma\in\Sigma}}{\Lambda_{\gamma}}))\leq k$ and so in particular, card$(g(\mathcal H_{\alpha}{_,}{_{\beta}})\Delta\mathcal H_{\alpha}{_,}{_{\beta}})\leq k$. But $\mathcal N_{_{\infty}}{^{^{0}}}$ is $k$-additive and $\{x\}\in\mathcal N_{_{\infty}}{^{^{0}}}$, so from above $g({\displaystyle\bigcup_{\gamma\in\Sigma}}{\Lambda_{\gamma}})\Delta\ ({\displaystyle\bigcup_{\gamma\in\Sigma}}{\Lambda_{\gamma}})\in \mathcal N_{_{\infty}}{^{^{0}}}$. Consequently, $g(\mathcal H_{\alpha}{_,}{_{\beta}})\Delta\mathcal H_{\alpha}{_,}{_\beta}\in\mathcal N_{_{\infty}}{^{^{0}}}$. Again as $\mathcal S_{0}\supseteq\mathcal N_{_{}\alpha}{^{^{0}}}$ $(\alpha<k)$ so $\mathcal H_{\alpha_{0}}{_,}{_{\beta_{0}}}\notin \mathcal N_{_{\infty}}{^{^{0}}}$; hence by Proposition $2.6$, $\mathcal H_{\alpha_{0}}{_,}{_{\beta_{0}}}$ is $(\mathcal S_{0}, \{\mathcal N_{_{\alpha}}{^{^{0}}}\}_{_{\alpha<k}})$-thick in $X$.\\
\vspace{.001cm}

We now proceed to define the class $\mathcal S_{1}$. We choose  $\mathcal H_{i}\in \{\mathcal H_{\alpha}{_,}{_{\beta}}: (\alpha,\beta)\neq(\alpha_{0},\beta_{0})\}$ $(i<\lambda)$ where $\lambda<k$ and for any $f\in\{0,1\}^{^{\lambda}}$ set $\mathcal H_{i}^{f(i)}= (X\setminus \mathcal H_{\alpha_{0}}{_,}{_{\beta_{0}}})\cap\mathcal H_{i}$ if $f(i)=0$ and $\mathcal H_{i}^{f(i)} = (X\setminus\mathcal H_{\alpha_{0}}{_,}{_{\beta_{0}}})\setminus\mathcal H_{i}$ if $f(i)=1$ and $\mathcal H_{f}=\displaystyle\bigcap_{i<\lambda}{\mathcal H_{i}^{f(i)}}$. Then the sets $\mathcal H_{f}$ for those $f\in \{0,1\}^{^{\lambda}}$ for which $\mathcal H_{f}\neq\emptyset$ along with $\mathcal H_{\alpha_{0}}{_,}{_{\beta_{0}}}$ constitutes a partition of $X$. We define $\mathcal S_{1}$ as the class of all sets which are of the form $(E\cap\mathcal H_{\alpha_{0}}{_,}{_{\beta_{0}}})\bigcup\displaystyle\bigcup_{f\in\{0,1\}^{^{\lambda}}}{(E_{f}\cap\mathcal H_{f})}$ where $E, E_{f}\in \mathcal S_{0}$, and $\mathcal H_{f}\neq\emptyset$. The Definition of $\mathcal S_{1}$ suggests that $\mathcal S_{1}\supseteq\mathcal S_{0}$ and so $\mathcal S_{1}$ is diffused. Now the complement of $(E\cap\mathcal H_{\alpha_{0}}{_,}{_{\beta_{0}}})\bigcup\displaystyle\bigcup_{f\in\{0,1\}^{^{\lambda}}}{(E_{f}\cap\mathcal H_{f})}$ is the set $\{(X\setminus E)\cap\mathcal H_{\alpha_{0}}{_,}{_{\beta_{0}}}\}\bigcup\displaystyle\bigcup_{f\in\{0,1\}^{^{\lambda}}}\{{(X\setminus E_{f})\cap\mathcal H_{f}}\}$ which again belongs to $\mathcal S_{1}$. Moreover, if $\{\lambda_{\alpha}:\alpha<k\}$ be a collection of ordinals where $\gamma, \lambda_{\alpha}<k$, then $\mu=\displaystyle\sum_{\alpha<\gamma}{\lambda_{\alpha}<k}$ because $k$ is regular and it is not hard to check that the union of the sets $(E^{\alpha}\cap\mathcal H_{\alpha_{0}}{_,}{_{\beta_{0}}})\bigcup\displaystyle\bigcup_{f\in\{0,1\}^{\lambda_{\alpha}}}{(E_{f}^{\alpha}\cap\mathcal H_{f}^{\alpha})}$ is a set of the form $((\displaystyle\bigcup_{\alpha<\gamma}{E^{\alpha}})\cap\mathcal H_{\alpha_{0}}{_,}{_{\beta_{0}}})\bigcup\displaystyle\bigcup_{\psi\in\{0,1\}^{^{\mu}}}(\mathcal H_{\psi}\cap E_{\psi})$, where $\mathcal H_{\psi}= {\displaystyle\bigcap_{i<\mu}{\mathcal H_{i}}^{\psi(i)}}$ and we take sets $\mathcal H_{\psi}$ for those $\psi\in \{0,1\}^{^{\mu}}$ for which $\mathcal H_{\psi}\neq \emptyset$. Thus $\mathcal S_{1}$ is a weakly $k$-additive algebra on $(X,G)$. Also from the calculations done in previous paragraph and the regularity of $k$, it follows that $\mathcal S_{1}$ is $G$-invariant.\\
\vspace{.1cm}

Now define $\mathcal N_{_{\alpha}}{^{^{1}}}=\{E\subseteq X : the \hspace{.08cm} (\mathcal S_{0},\{\mathcal N_{_{\alpha}}{^{^{0}}}\}_{\alpha<k})-cover\hspace{.01cm} of\hspace{.1cm} {E\cap\mathcal H_{\alpha_{0}}{_,}{_{\beta_{0}}}}\hspace{.08cm} is \hspace{.2cm} a\hspace{.1cm} member\hspace{.1cm} of\hspace{.1cm}\\ \mathcal N_{_{\alpha}}{^{^{0}}}\}$. From the definition of $(\mathcal S_{0},\{\mathcal N_{_{\alpha}}{^{^{0}}}\})$ cover and conditions (iii) and (iv) of Definition $2.2$ it follows that $\mathcal N_{_{\alpha}}{^{^{0}}}\subseteq\mathcal N_{_{\alpha}}{^{^{1}}}$. The same observations establish conditions (i), (iii), (iv), (v) and (vi) of Definition $2.2$ for $\mathcal N_{_{\alpha}}{^{^{1}}}$. To check that $\{\mathcal N_{_{\alpha}}{^{^{1}}}\}_{_{\alpha<k}}$ is a $k$-small system on $(X,G)$, we need only verify condition (ii) of Definition $2.2$. Let $E\in\mathcal N_{_{\alpha}}{^{^{1}}}$ and $\mathcal C$ be the $(\mathcal S_{0},\{\mathcal N_{_{\alpha}}{^{^{0}}}\}_{\alpha<k})$-cover of $E\cap\mathcal H_{\alpha_{0}}{_,}{_{\beta_{0}}}$. Then $\mathcal C\in \mathcal N_{_{\alpha}}{^{^{0}}}$ and therefore $g{\mathcal C}\in\mathcal N_{_{\alpha}}{^{^{0}}}$ by condition (ii) of Definition $2.2$. But $gE\cap\mathcal H_{\alpha_{0}}{_,}{_{\beta_{0}}}\subseteq(g(\mathcal H_{\alpha_{0}}{_,}{_{\beta_{0}}})\Delta\mathcal H_{\alpha_{0}}{_,}{_{\beta_{0}}})\Delta\hspace{.01cm} g(E\cap\mathcal H_{\alpha_{0}}{_,}{_{\beta_{0}}})$ where $g(\mathcal H_{\alpha_{0}}{_,}{_{\beta_{0}}})\Delta\mathcal H_{\alpha_{0}}{_,}{_{\beta_{0}}}\in\mathcal N_{_{\infty}}{^{^{0}}}$ and $g(E\cap\mathcal H_{\alpha_{0}}{_,}{_{\beta_{0}}})\subseteq g{\mathcal C}\in \mathcal N_{_{\alpha}}{^{^{0}}}$. Hence the $(\mathcal S_{0},\{\mathcal N_{_{\alpha}}{^{^{0}}}\}_{\alpha<k})$-cover of $gE\cap\mathcal H_{\alpha_{0}}{_,}{_{\beta_{0}}}$ is a member of $\mathcal N_{_{\alpha}}{^{^{0}}}$ by the same reasoning as above.\
\vspace{.1cm}

The small system $\{\mathcal N_{_{\alpha}}{^{^{1}}}\}_{\alpha<k}$ is upper semicontinuous relative to $\mathcal S_{1}$. Let $\{\Gamma^{^{\alpha}}\}_{\alpha<k}$ be a nested $k$-sequence of sets from $\mathcal S_{1}$ such that for some $\alpha_{_{0}}<k$, $\Gamma^{^{\alpha}}\notin\mathcal N_{_{\alpha{_{0}}}}{^{1}}$ for all $\alpha<k$. Let $\mathcal C^{^{\alpha}}$ be the $(\mathcal S_{0},\{\mathcal N_{_{\alpha}}{^{^{0}}}\}_{\alpha<k})$-cover of $\Gamma^{^{\alpha}}\cap\mathcal H_{\alpha_{0}}{_,}{_{\beta_{0}}}$. We write  $\Gamma^{^{\alpha}}= (E^{\alpha}\cap\mathcal H_{\alpha_{0}}{_,}{_{\beta_{0}}})\bigcup\displaystyle\bigcup_{f\in\{0,1\}^{\lambda_{\alpha}}}(E_{f}^{\alpha}\cap\mathcal H_{f}^{\alpha})$. Since $\Gamma^{^{\alpha}}\cap\mathcal H_{\alpha_{0}}{_,}{_{\beta_{0}}}= E^{\alpha}\cap\mathcal H_{\alpha_{0}}{_,}{_{\beta_{0}}}$, $\mathcal C^{^{\alpha}}$ is also the $(\mathcal S_{0},\{\mathcal N_{_{\alpha}}{^{^{0}}}\}_{\alpha<k})$-cover of $E^{\alpha}\cap\mathcal H_{\alpha_{0}}{_,}{_{\beta_{0}}}$. Again, $\mathcal H_{\alpha_{0}}{_,}{_{\beta_{0}}}$ is $(\mathcal S_{0},\{\mathcal N_{_{\alpha}}{^{^{0}}}\}_{\alpha<k})$-thick in $X$, so $E^{\alpha}$ is identical with $\mathcal C^{^{\alpha}}$. But $\mathcal C^{^{\alpha}}\notin \mathcal N_{_{\alpha{_{0}}}}{^{0}}$, so $E^{{\alpha}}\notin \mathcal N_{_{\alpha{_{0}}}}{^{0}}$. We set $F^{\alpha}= \displaystyle\bigcap_{\beta<\alpha}{E^{\beta}}$. Then $E^{\alpha}\setminus F^{\alpha}\in \mathcal N_{_{\infty}}{^{^{0}}}$. So $F^{\alpha}\notin \mathcal N_{_{\alpha{_{0}}}}{^{0}}$. Also $F^{\alpha}\in \mathcal S_{0}$ $(\alpha<k)$ and $\{F^{\alpha}\}_{_{\alpha<k}}$ is a nested $k$-sequence. Therefore $\displaystyle\bigcap_{\alpha<k}{F^{\alpha}}\notin\mathcal N_{_{\infty}}{^{^{0}}}$ because $\{\mathcal N_{_{\alpha}}{^{^{0}}}\}_{\alpha<k}$ is upper semicontinuous relative to $\mathcal S_{0}$. Hence $\displaystyle\bigcap_{\alpha<k}{E^{\alpha}}\notin\mathcal N_{_{\infty}}{^{^{0}}}$. But $\displaystyle\bigcap_{\alpha<k}{E^{\alpha}}$ is the $(\mathcal S_{0},\{\mathcal N_{_{\alpha}}{^{^{0}}}\}_{\alpha<k})$-cover of $(\displaystyle\bigcap_{\alpha<k}{E^{\alpha}})\cap\mathcal H_{\alpha_{0}}{_,}{_{\beta_{0}}}$ and hence the $(\mathcal S_{0},\{\mathcal N_{_{\alpha}}{^{^{0}}}\}_{\alpha<k})$ cover of $(\displaystyle\bigcap_{\alpha<k}{\Gamma^{\alpha}})\cap\mathcal H_{\alpha_{0}}{_,}{_{\beta_{0}}}$. Therefore $\displaystyle\bigcap_{\alpha<k}{\Gamma^{\alpha}}\notin \mathcal N_{_{\infty}}{^{^{1}}} $.\
\vspace{.1cm}

We  now check that the $k$-additive algebra $\mathcal S_{1}$ is admissible with respect to $\{\mathcal N_{_{\alpha}}{^{^{1}}}\}_{\alpha<k}$. Since $\mathcal S_{0}$ is admissible with respect to $\{\mathcal N_{_{\alpha}}{^{^{0}}}\}_{\alpha<k}$, by condition (iii) of Definition $2.2$, it follows that $\mathcal H_{\alpha_{0}}{_,}{_{\beta_{0}}}\in\mathcal S_{1}\setminus\mathcal N_{_{\alpha}}{^{^{1}}}$ for some $\alpha<k$. Again, $(X\setminus\mathcal H_{\alpha_{0}}{_,}{_{\beta_{0}}})\cap\mathcal H_{\alpha}{_,}{_{\beta}}\in \mathcal S_{1}\cap\mathcal N_{_{\alpha}}{^{^{1}}}$ $(\alpha<k)$ for every $\alpha<k$ provided $(\alpha,\beta)\neq(\alpha_{_{0}},\beta_{_{0}})$. So condition (i) of Definition $2.3$ is verified. Now let $M\in \mathcal N_{_{\alpha}}{^{^{1}}}$ and $N$ be the $(\mathcal S_{0},\{\mathcal N_{_{\alpha}}{^{^{0}}}\}_{\alpha<k})$-cover of $M\cap\mathcal H_{\alpha_{0}}{_,}{_{\beta_{0}}}$. Then $N \in {\mathcal S_{0}} \cap {\mathcal N_{_{\alpha}}{^{^{0}}}}$ and we set $P= N\cup (X\setminus\mathcal H_{\alpha_{0}}{_,}{_{\beta_{0}}})$. Evidently, $M\subseteq P\in\mathcal S_{1}\cap\mathcal N_{_{\alpha}}{^{^{1}}}$. This verifies condition (ii) of Definition $2.3$. Finally, let us choose an arbitrary collection $\{E_{\beta}: \beta\in\Lambda\}$ of mutually disjoint sets from $\mathcal S_{1}\setminus\mathcal N_{_{\alpha}}{^{^{1}}}$. Then the sets $E_{\beta}\cap\mathcal H_{\alpha_{0}}{_,}{_{\beta_{0}}}$ are also mutually disjoint and for each $\beta\in\Lambda$, there exists $F_{\beta}$ where $F_{\beta}\in\mathcal S_{0}$ such that $E_{\beta}\cap\mathcal H_{\alpha_{0}}{_,}{_{\beta_{0}}}= F_{\beta}\cap\mathcal H_{\alpha_{0}}{_,}{_{\beta_{0}}}$. Now define sets $G_{\beta}= F_{\beta}\setminus\displaystyle\bigcup_{\gamma<\beta}{F_{\gamma}}$ $(\beta,\gamma\in\Lambda)$, which are certainly members of $\mathcal S_{0}$. Let $\mathcal C_{\beta}$ be the $(\mathcal S_{0},\{\mathcal N_{_{\alpha}}{^{^{0}}}\}_{\alpha<k})$-cover of $E_{\beta}\cap\mathcal H_{\alpha_{0}}{_,}{_{\beta_{0}}}$. Then it is also the $(\mathcal S_{0},\{\mathcal N_{_{\alpha}}{^{^{0}}}\}_{\alpha<k})$-cover of $G_{\beta}\cap\mathcal H_{\alpha_{0}}{_,}{_{\beta_{0}}}$ and so $C_{\beta}\setminus G_{\beta}\in\mathcal N_{_{\infty}}{^{^{0}}}$. We set $D_{\beta}=\mathcal C_{\beta}\setminus(\mathcal C_{\beta}\setminus G_{\beta})$. Then $E_{\beta}\cap\mathcal H_{\alpha_{0}}{_,}{_{\beta_{0}}}\subseteq D_{\beta}\subseteq G_{\beta}$, $D_{\beta}\in \mathcal S_{0}\setminus \mathcal N_{_{\alpha}}{^{^{0}}}$ and $D_{\beta}$ $(\beta\in\Lambda)$ are mutually disjont. By admissibility of $\mathcal S_{0}$ with respect to $\{\mathcal N_{_{\alpha}}{^{^{0}}}\}_{\alpha<k}$, Card($\Lambda$)$\leq k$ which verifies condition (iii) of Definition $2.3$.\
\vspace{.5cm}

In a similar manner, we define a class $\mathcal S_{2}$ as consisting of those sets which have the form $(E\cap(X\setminus\mathcal H_{\alpha_{0}}{_,}{_{\beta_{0}}}))\bigcup\displaystyle\bigcup_{f\in\{0,1\}^{\lambda}}(E_{f}\cap\mathcal H_{f})$ where $E,E_{f}\in \mathcal S_{0}$, $\lambda<k$, $\mathcal H_{f}= \bigcap{\mathcal H_{_{i}}{^{f(i)}}}$ and $\mathcal H_{_{i}}{^{f(i)}}= \mathcal H_{\alpha_{0}}{_,}{_{\beta_{0}}}\cap\mathcal H_{_{i}}$ if $f(i)=0$ and $\mathcal H_{_{i}}{^{f(i)}}= \mathcal H_{\alpha_{0}}{_,}{_{\beta_{0}}}\setminus\mathcal H_{_{i}}$ if $f(i)=1$ where $\mathcal H_{_{i}}\in\{\mathcal H_{\alpha_{0}}{_,}{_{\beta_{0}}}\cap\mathcal H_{\alpha}{_,}{_{\beta}} : (\alpha,\beta)\neq (\alpha_{_{0}},\beta_{_{0}})\}$; the union being taken over all $f\in\{0,1\}^{\lambda}$ for which $\mathcal H_{f}\neq\emptyset$. Using a procedure similar as above, it can be verified that $\mathcal S_{2}$ is a weakly $k$-additive algebra on $(X,G)$ which is diffused because it contains $\mathcal S_{0}$. It is evident from the constructions of $\mathcal S_{1}$ and $\mathcal S_{2}$ that they both are $k$-generated if $\mathcal S_{0}$ is so. Moreover, the class $\mathcal N_{_{\alpha}}{^{^{2}}}=\{E\subseteq X : the \hspace{.06cm} (\mathcal S_{0},\{\mathcal N_{_{\alpha}}{^{^{0}}}\}_{\alpha<k})-cover\hspace{.06cm} of\hspace{.1cm} {E\cap(X\setminus\mathcal H_{\alpha_{0}}{_,}{_{\beta_{0}}})}\hspace{.06cm} is \hspace{.06cm} a\hspace{.1cm} member\hspace{.1cm} of\hspace{.1cm}\\ \mathcal N_{_{\alpha}}{^{^{0}}}\}$ is a $k$-small system on $(X,G)$ such that $\mathcal S_{2}$ is admissible with respect to $\{\mathcal N_{_{\alpha}}{^{^{2}}}\}_{\alpha<k}$ and $\{\mathcal N_{_{\alpha}}{^{^{2}}}\}_{\alpha<k}$ is upper semicontinuous relative to $\mathcal S_{2}$. But there can be no $k$-additive algebra $\mathcal S$ on $(X,G)$ containing $\mathcal S_{1}\cup\mathcal S_{2}$ which is admissible with respect to any upper semicontinuous $k$-small system on $(X,G)$ because if so, then it would contain the class $\{\mathcal H_{_{\alpha,\beta}} : \alpha<k, \beta<k\}$. But this would violate Theorem $2.7$.\\
\vspace{.02cm}

\textbf{NOTE :} Kharazishvili gave another generalization $[8]$ of Pelc and Prikry's theorem based on the fact that under Martin's axiom there exists absolutely nonmeasurable functions.\\
\vspace{.08cm}
\hspace{1cm}

\begin{center}

\end{center}
\vspace{.1cm}


{\normalsize
\textbf{S.Basu}\\
\vspace{.02cm}
\hspace{.35cm}
\textbf{Dept of Mathematics}\\
\vspace{.02cm}
\hspace{.35cm}
\textbf{Bethune College , Kolkata} \\
\vspace{.01cm}
\hspace{.35cm}
\textbf{W.B. India}\\
\vspace{.02cm}
\hspace{.1cm}
\hspace{.29cm}\textbf{{e-mail : sanjibbasu08@gmail.com}}\\

\textbf{D.Sen}\\
\vspace{.02cm}
\hspace{.35cm}
\textbf{Saptagram Adarsha vidyapith (High), Habra , $24$ Parganas (N)} \\
\vspace{.01cm}
\hspace{.35cm}
\textbf{W.B. India}\\
\vspace{.02cm}
\hspace{.1cm}
\hspace{.29cm}\textbf{{e-mail : reachtodebasish@gmail.com}}
}

\end{document}